\documentclass[12pt,a4paper,twoside,final,notitlepage, leqno]{article}
\usepackage[english]{babel}
\usepackage[T1]{fontenc}  
\usepackage{graphicx}
\usepackage{amsmath,amsthm,epsfig,amsfonts,bbm}  

\newcommand{\pequationdeb}{$$ \left\{ \begin{minipage}[c]{130mm}}
\newcommand{\pequationfin}{\end{minipage}
                           \right. $$}

\def \smb {{\scriptstyle \bullet }}
\newcommand{\monitem}{ \smallskip \noindent $\bullet$ \quad  } 
\newcommand{\moneq}{\vspace*{-6pt} \begin{equation} \displaystyle } 
\newcommand{\moneqstar}{\vspace*{-6pt} \begin{equation*} \displaystyle } 
\newcommand{\monendstar}{\vspace*{-6pt} \end{equation*}   }
\newcommand{\monend}{\vspace*{-6pt} \end{equation}   }
\newcommand{\beq}     {\begin{equation}}
\newcommand{\enq}     {\end{equation}}
\newcommand{\be}    {\begin{enumerate}}
\newcommand{\ee}    {\end{enumerate}}

\newcommand{\Bb}

\def\br {\break}

\def\section*#1{}
\def\resume{\if@twocolumn
\section*{R\'esum\'e}
\else \small
\quotation{\bf \it R\'esum\'e \rule[1mm]{1.5mm}{0.2mm}\vspace{0pt}}
\fi}
\def\endresume{\if@twocolumn\else\endquotation\fi}
\def\abstract{\if@twocolumn
\noindent\section*{{\bf Abstract}}
\else \small
\quotation{\noindent \bf {Abstract.} \rule[1mm]{1.5mm}{0.2mm}\vspace{0pt}}
\fi}
\def\endabstract{\if@twocolumn\else\endquotation\fi}

\usepackage{fancyhdr}
\renewcommand{\headrulewidth}{0pt}
\begin{document}

\fancypagestyle{plain}{ \fancyfoot{} \renewcommand{\footrulewidth}{0pt}}
\fancypagestyle{plain}{ \fancyhead{} \renewcommand{\headrulewidth}{0pt}} 

~
 
~

\bigskip \bigskip   \bigskip

\centerline {\bf \LARGE   Comparison of Simulations of Convective Flows} 

\bigskip \bigskip \bigskip

\centerline { \large  Pierre Lallemand$^{a}$ and Fran\c{c}ois Dubois$^{b c}$}   

\bigskip  

\centerline { \it  \small    $^a$   Beijing Computational Science Research Center,   } 
\centerline { \it  \small      Beijing Run Ze Jia Ye, China.   } 
\centerline { \it  \small    $^b$   Conservatoire National des Arts et M\'etiers,  Paris, France,  } 
\centerline { \it  \small   Laboratoire de M\'ecanique des Structures 
et des Syst\`emes Coupl\'es.  } 
\centerline { \it  \small  $^c$   Department of Mathematics, University  Paris-Sud,}  
\centerline { \it  \small  B\^at. 425, F-91405 Orsay Cedex, France.} 
\centerline { \it  \small pierre.lallemand1@free.fr, francois.dubois@math.u-psud.fr}    \bigskip   
 
\bigskip


\centerline {  15 may 2015~\footnote{~Contribution published in   
  {\it  Communications in Computational Physics}, doi: 10.4208/cicp.2014.m400, june 2015. 
presented at the 10th  International Conference for Mesoscopic Methods in Engineering and Science,
 Oxford, UK,   22-26 July 2013. }} 
 
\bigskip

\bigskip 
\noindent  {\bf Abstract. } \qquad 
We show   that a  single  particle distribution for the     
``energy-conserving''  
D2Q13 lattice Boltzmann scheme can simulate coupled 
effects involving advection and diffusion of velocity and temperature. 
We consider various test cases:  non-linear waves with periodic boundary conditions, 
a test case with buoyancy, propagation of transverse waves, Couette and Poiseuille flows.
We test various boundary conditions and propose to mix bounce-back and 
anti-bounce-back numerical boundary conditions to take into account
velocity and temperature Dirichlet conditions. 
We present also  first results for the  de Vahl Davis heated cavity.
Our results are compared with  the coupled 
D2Q9-D2Q5 lattice Boltzmann approach for the Boussinesq system 
and with an elementary finite differences solver for the compressible Navier-Stokes
equations. Our main experimental result is the loss of symmetry in the de Vahl Davis cavity 
computed with the single D2Q13 lattice Boltzmann model without the Boussinesq hypothesis.
This result is confirmed by a direct Navier Stokes simulation with finite differences.
 
 $ $ \\  [2mm]

%

\bigskip \bigskip  \newpage \noindent {\bf \large   Introduction}  
 
\fancyfoot[C]{\oldstylenums{\thepage}}
\fancyhead[OC]{\sc{Comparison of Simulations of Convective Flows }}

%

\monitem 
Lattice Boltzmann schemes have proven their efficiency for the computation of 
quasi-incompressible flows. We refer {\it e.g.} to \cite{BSV92, JKL05, LL00} among others. 
In these cases, the physical conservations of  mass and momentum are implemented
in the framework of lattice Boltzmann schemes. 
When compressible effects are taken into account, it is necessary to add the  conservation
of energy. A classical approach is to begin with weakly compressible effects that can be 
modelled with the so-called Boussinesq approximation. In this case, the incompressibility condition 
remains a good approximation and coupled effects between conservations of
momentum and energy are taken into account with a precise thermodynamical analysis.
We refer to Landau \cite{LL59}  or Batchelor \cite{Ba67}  
for the derivation of the Boussinesq approximation. 
The implementation of the Boussinesq approximation is possible with the lattice Boltzmann
approach with the introdution of two particle distributions. 
This idea has been also proposed in the context of finite volumes by the 
team of Perthame \cite{KP94}, 
and with lattice Boltzmann schemes by 
Eggels and Somers  \cite {ES95}, 
Mezrhab {\it et al} \cite {MMJN10} and    
  Wang {\it et al} \cite {WWLL13}  
among others. 

\monitem 
In this contribution, we study a direct approximation of the compressible Navier Stokes equations
with an ``energy-conserving''  
lattice Boltzmann scheme using a single particle distribution.   
A first tentative study  \cite{LL03} has shown that 
for a critical value of the Prandlt number, 
the thermal wave and the viscous one merge together, the physics is badly represented
and an instability occurs in general. 
In consequence, no satisfying compressible flows have been obtained with this direct  numerical modelling.
In a second tentative  \cite{LD13}, we have analyzed with great details 
several lattice Boltzmann schemes 
with  four conservation laws in two space dimensions. 
With an adequate fitting of the parameters of the scheme, it is possible to enlarge the zone   
in the spectral space where the thermal and viscous waves remain decoupled.  
Moreover, these parameters  guarantee also the isotropy                            
of the acoustic waves. Our objective is to enlarge the domain of validity of our previous study: 
incorporate the treatment of boundary conditions with rigid walls with a given temperature    
or adiabatic boundaries, study several couplings between velocity and temperature
for elementary Couette and Poiseuille flows, study the possibility of Dirichlet 
and Neumann boundary conditions. Finally, our objective is the simulation 
of the  de Vahl Davis  test case \cite{dVD83}   described in Figure~1.

\monitem 
The outlook of the article is the following. 
In Section~1, we recall fundamental aspects of the  coupled 
D2Q9-D2Q5 lattice Boltzmann approach. We present our actual choices for
the implementation of the lattice Boltzmann  approach with the 
D2Q13 stencil and to  treat all the physical fields 
with  single particle distribution and the    
D2Q13 scheme.  
In Section~3, we develop 
a very elementary finite-difference approach  for  the compressible    
Navier-Stokes equations. With this tool, we can  
compare our new D2Q13 approach with a classical reference. 
In Section~4, we consider a simple test case for non-linear waves.
We study the  buoyancy in Section~5, the propagation of transverse waves in Section~6, 
the simulation of Couette flows in Section~7
and Poiseuille flows in Section~8.  
In Section~9, we consider a test case to take into account various
temperature and flux-type boundary conditions. 
First results for the  de Vahl Davis heated cavity are presented in Section~10. 

\smallskip \smallskip 
\centerline {\includegraphics[width=.50 \textwidth]   {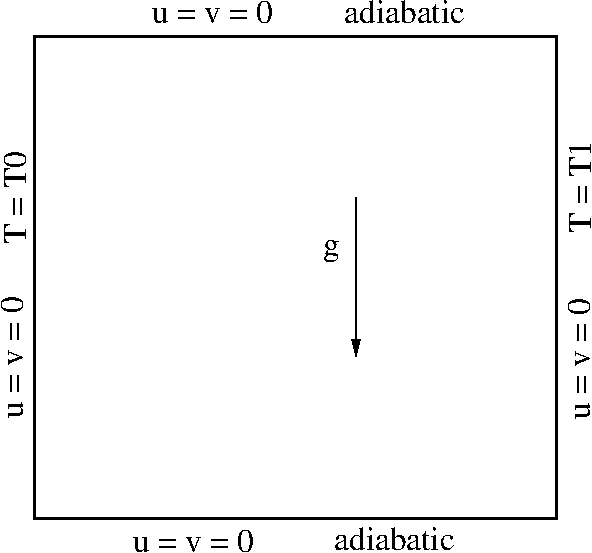} }  

\smallskip \noindent  {\bf Figure 1}. \quad 
De Vahl Davis test case for natural convection 
\smallskip  \smallskip 

\bigskip \bigskip  \noindent {\bf \large 1) \quad Coupled D2Q9-D2Q5 lattice Boltzmann scheme}  

\monitem 
The  Boussinesq approximation  of the compressible Navier-Stokes equations 
can be written as a system of coupled partial differential equations. 
The unknowns are the vector field of velocity $ \, u ,\ $  
and the scalar fields of temperature $ \, T  \,$ and  pressure~$\, p$. 
The parameters are the shear viscosity $\, \nu$, the temperature dissipation rate   
$ \, \kappa , \,$ the thermal expansion coefficient $ \, \beta $ and $ \, g  \,$ 
the Earth's gravity. The buoyancy term $ \, ( 1 - \beta  \, (T - T_0)  ) \, g \,$ is a source term
for the momentum equation and the velocity field $ \, u \, $ directly imposes strong
constraints for the transport of temperature. Assuming that the density is $\rho=1$,            
the equations of the Boussinesq system are 
\moneq  \label{intro-boussinesq} 
 \left\{ \begin{array}{rl}  \displaystyle 
{\rm div}  u & =\,  0 \,, \\ \displaystyle 
 {{\partial u}\over{\partial t}}  + u  \smb  \nabla u + \nabla p 
- \nu \triangle u & =  \, \big( 1 - \beta  \, (T - T_0) \big) \, g  \,, \\  \displaystyle 
 {{\partial T}\over{\partial t}}  + u \smb \nabla T  
- \kappa \triangle T  & = \, 0  \, . 
\end{array}  \right. \monend
The    Rayleigh number is defined from the temperature difference 
$ \, \Delta T \,\equiv \, T_1 - T_0 \,$  between the two sides according to 
\moneq  \label{intro-rayleigh} 
 {\rm  R_a}  \equiv {{ \mid  g   \mid  \, \beta \, \Delta T \, L^3}
 \over{\nu \, \kappa}}  \, . 
\monend 
%
%

\monitem 
A difficult stationary test case is the computation 
of the velocity and temperature fields for $ \,  {\rm  R_a} = 10^6 $.    
The references are the original contribution 
of de Vahl Davis \cite{dVD83}, the     
Le Qu\'er\'e  \cite{LQ91} and the associated workshop in the 2000's  
with various  Navier-Stokes  solvers, 
the introduction of the  D2Q9-D2Q5  coupled approximation 
by Mezrhab {\it et al.} \cite{MMJN10}
and the very precise results of 
Wang {\it et al}  \cite{WWLL13}  with the same approach.


\monitem  
Recall that the discrete velocities of a D2Q5 lattice Boltzmann scheme
follow the axis of coordinates:
\moneq  \label{d2q5-vitesses} 
v_j \, \in \, \{ (0,\ 0) , \, (1,\ 0) , \,  (0, \ 1) , \, 
 (-1,\ 0) , \,  (0, \ -1) \} \,, \quad 0 \leq j \leq 4 \, . 
\monend 
For a D2Q9 scheme we add to the previous D2Q5 velocities (\ref{d2q5-vitesses}) 
the four ones along the diagonals: 
\moneq  \label{d2q9-vitesses} 
v_j \, \in \, \{ (1,\ 1) , \,  (-1,\ 1) , \, 
 (-1,\ -1) , \,  (1,\ -1) \} \,, \quad 5 \leq j \leq 8 \, . 
\monend 
The flow is simulated with a  D2Q9 lattice Boltzmann scheme 
with 3 conserved moments, 
the density and the two components of the momentum: 
\moneq  \label{d2q9-conserve} 
\rho \equiv \sum_{j=0}^8 f_j \,, \quad  
 (j_x \,,\, j_y )  \equiv \sum_{j=0}^8  v_j \, f_j \,. 
\monend 
The six other moments of the fluid are presented in the reference  \cite{LL00}. 
The equilibrium values for the  moments of order two  
have to take into account the compressible effects:
\moneq  \label{d2q9-equilibres} 
E^{\rm eq}  \,= \,  \alpha \rho + 3\ \frac{j_x^2+j_y^2}{\rho} \,, \quad 
XX^{\rm eq} \,= \,   \frac{j_x^2-j_y^2}{ \rho} \,, \quad 
XY^{\rm eq}  \,= \,   \frac{j_x j_y}{ \rho} \, . 
\monend 
The equilibrium properties either have no influence on the physical properties or are set          
to give an isotropic shear viscosity.                                                             
The sound velocity $\, c_s$, the shear viscosity $ \, \mu \,$
and the bulk viscosity $ \, \zeta \,$ are given from the previous equilibria 
according to
\moneq  \label{d2q9-physic}  
c_s  \,= \,   \sqrt{{4+\alpha}\over{6}} \,, \quad 
\mu  \,= \, \frac{1}{3} \, \Big( \frac{1}{s_{XX}}-\frac{1}{2}  \Big)  \,, \quad 
\zeta  \,= \, -\alpha  \,  \Big( \frac{1}{s_{E}}-\frac{1}{2} \Big) \, . 
\monend

\monitem  
The temperature is simulated with a simple D2Q5 scheme 
with only one  conserved moment
\moneq  \label{d2q5-conserve} 
 T \equiv \sum_{j=0}^4 g_j \, . 
\monend 
The other nontrivial equilibrium values follow the relations
\moneq  \label{d2q5-equil} 
E^{\rm eq}   \,= \,   \beta \, \rho \,, \quad 
j_x^{\rm eq} \,= \,   \rho \,  V_x \,, \quad 
j_y^{\rm eq}  \,= \,  \rho \,  V_y  \,, \quad 
XX^{\rm eq}  \,= \,  0  \, . 
\monend
The diffusion coefficient $ \, \kappa \,$ is easy to identify: 
\moneq  \label{d2q5-kappa} 
\kappa   \,= \,   \frac{\beta + 4}{10} \Big( \frac{1}{s_E}-\frac{1}{2} \Big)  \, . 
\monend
This D2Q5 model as defined does not satisfy the  Galilean invariance with respect to advection  
at uniform speed $\{V_x,\ V_y\}$ (see Qian and Zhou \cite{QZ98}); 
the equivalent equation  for the D2Q5 scalar scheme is equal to 
\moneqstar 
{{\partial T}\over{\partial t}} \,+\, V_x \,  {{\partial T}\over{\partial x}} 
\, + \, V_y \,  {{\partial T}\over{\partial y}}  
- \kappa \, \triangle T \, + 
\Big( \frac{1}{s_E}-\frac{1}{2} \Big) \Big( V_x^2 \, {{\partial^2 T}\over{\partial x^2}} 
\,+\, 2 \, V_x \, V_y \, {{\partial^2 T}\over{\partial x  \, \partial y}}  
\, + \, V_y^2 \, {{\partial^2 T}\over{\partial y^2}} \Big)  = 0 
\monendstar 
and can be  easily identified with the methods
developed {\it e.g.} in \cite{DL09}.

\bigskip \bigskip  \newpage 
\noindent {\bf \large 2) \quad Compressible D2Q13 lattice Boltzmann scheme}  

\monitem 
The stencil of the D2Q13 lattice Boltzmann scheme 
is built (see {\it e.g.} \cite{LL03,LD13})
on the D2Q9 scheme with the following complementary velocity set:
\moneq  \label{d2q13-vitesses} 
v_j \, \in \, \{ (2 ,\ 0) , \,  (0 ,\ 2 ) , \, 
 (-2 ,\ 0) , \,  (0 ,\ -2 ) \} \,, \quad 9 \leq j \leq 12 \, . 
\monend 
%
A family of 13 orthogonal moments are generated by an elementary linear
mapping of the particle distribution $ \, f_j $:
\moneqstar  
m_k \,=\, \sum_{j=0}^{12} M_{kj} \, f_j  \, . 
\monendstar 
The coefficients $ \,  M_{kj} \, $ of the matrix are computed 
from the 13 velocities presented in  (\ref{d2q5-vitesses}),  (\ref{d2q9-vitesses})
and (\ref{d2q13-vitesses}) with the help of 
polynomials $ \, p_k \,$ by the condition
\moneq  \label{d2q13-moments} 
 M_{kj} \, =\, p_k (v^x_j, \, v^y_j) \,, \qquad 0 \leq j, \, k \leq 8 \,.
\monend
%
The following set $ \, \{ p_k \} \,$  of polynomials 
are presented in   (\ref{d2q13-polynomes}) as combinations of monomials of increasing power. 
They have  been chosen as symmetric as possible and have been 
orthogonalized. Instead of giving the final moment matrix, we give the ``recipe''
to build it in terms of the components $x \equiv v^x_j$ and $y \equiv v^y_j$ of the 13 basic velocities. 
%
\moneq  \label{d2q13-polynomes} 
 \left\{ \begin{tabular}{ccc}
scalars & $\rho$ & 1 \cr
       &  $E$   & $-28+13 \, (x^2+y^2)$    \cr
       &  $\epsilon$ & $140+(x^2+y^2) \, (-361/2+77(x^2+y^2)/2)$ \cr
  &  $\varpi$  & $-12+(x^2+y^2)(\frac{581}{12}+(x^2+y^2)(-\frac{273}{8}+\frac{137}{24}(x^2+y^2)))$  \cr
vectors & $j_x$ & $x$ \cr
       & $j_y$ & $y$ \cr
       & $q_x$ & $x \, (3+x^2+y^2)$ \cr
       & $q_y$ & $y \, (3+x^2+y^2)$ \cr
       & $r_x$ & $x \, (\frac{101}{6}+(x^2+y^2)(-\frac{63}{4}+\frac{35}{12}(x^2+y^2)))$ \cr
       & $r_y$ & $y \, (\frac{101}{6}+(x^2+y^2)(-\frac{63}{4}+\frac{35}{12}(x^2+y^2)))$ \cr
tensors & $XX$ & $x^2-y^2$  \cr
        & $XY$ & $ x\ y $  \cr
        & $XX_e$ & $(x^2-y^2) \, (-\frac{65}{12}+\frac{17}{12}(x^2+y^2))\, . $   \cr
\end{tabular}  \right. \monend

\monitem 
The collisions conserve two scalars $\rho$ and $E$ 
and two vector components $j_x$ and $j_y$.
The equilibrium values of the other moments 
and the relaxation rates are constrained by the result 
of a linearized analysis of the four hydrodynamic modes. 
The four  modes show isotropic behaviour for their attenuation 
and propagation velocity meaning that Galilean invariance
is achieved. 
The equilibrium expressions can be taken as simple functions of the conserved variables
that have the same symmetry properties.  The choice of linear and quadratic expressions leads to:     
%
\moneq  \label{d2q13-moments-equil} 
 \left\{   \begin{array}{rlrl} 
q_x^{\rm eq}   \,&= \,   j_x \, (c_1+h_1 \rho + k_1 E) \,, 
& r_x^{\rm eq}   \,&= \,   j_x \, (c_2+h_2 \rho + k_2 E) \,, \cr 
 \epsilon^{\rm eq} \,&= \, c_{\epsilon\rho} \, \rho + c_{\epsilon E} \, E \,,   
& \varpi^{\rm eq} \,&= \,   c_{\varpi \rho} \, \rho + c_{\varpi E} \,  E \,,  \cr 
XX^{\rm eq}  \,&= \,  \frac{j_x^2-j_y^2}{ \rho}\,,  
& XY^{\rm eq}   \,&= \, \frac{j_x j_y}{ \rho} \,,  \qquad 
 XX_e^{\rm eq}   \, = \, 0 \,  . 
\end{array}  \right. \monend
With this specific choice of $ \, \epsilon^{\rm eq} \,$ and 
$ \, \varpi^{\rm eq} \,$ we introduce new parameters 
that give some freedom to develop as in  \cite{LD13} a solution to 
the unphysical coupling observed in \cite{LL03}. 

\monitem 
Following  H\'enon  \cite{He87}, many formulae can be simplified using  
%
\moneq  \label{d2q13-henon} 
\sigma_i \equiv {{1}\over{s_i}} - {1\over2}  \,.   
\monend
With an asymptotic analysis as the one presented in \cite{LD13}, 
we can recover the physical waves : 
two acoustics, one transverse and  one longitudinal diffusion.
To satisfy correct advection of these four waves in the presence on a uniform  
background velocity, the following relationships have to be satisfied :        
\moneq  \label{d2q13-diffusion} 
 \left\{ \begin{array} {l}  \displaystyle 
h_1 = \frac{17}{26}-\frac{c_1}{2}-\frac{E_0}{13}  \,, \qquad 
k_1 = \frac{2}{13}   \,, \\  \vspace{-.4cm}    \\ \displaystyle 
h_2 = -\Big(\, \frac{39}{2}+\frac{13}{2}c_1+E_0 \, \Big)\ k_2   \\ \displaystyle 
\quad  \displaystyle 
-\frac{7}{624}(13\ c_1+95+2\ E_0) \, \frac{874481 + 459459\ c_1 -103428\ c_{\epsilon E}+70686\ E_0}
{114404+51051\ c_1 - 11492\ c_{\epsilon E} + 7854\ E_0} \,   
\end{array}  \right. \monend
for a situation with density equal to 1 and ``energy" equal to $E_0$.                 
In order to enforce isotropy at second order around a null velocity, we have  to set    
\moneq  \label{d2q13-sigma-qx} 
 \left\{ \begin{array} {l}  \displaystyle 
\sigma_{qx} = -\frac{1309}{2} \, \frac{\sigma_{XX} \,  (13\ c_1 + 95 + 2\ E_0)}
{114404 + 51051\ c_1 - 11492\ c_{\epsilon E} + 7854\ E_0}  \,, \\ \displaystyle 
c_{\epsilon \rho} = 140+28\ c_{\epsilon E} \,  \\ \displaystyle 
\qquad + \frac{(13\ c_1 + 95+2\ E_0) \, (114404+51051\ c_1 - 11492\ c_{\epsilon E} + 7854\ E_0)}
{22984 \ Pr}  \,,  \\ \displaystyle 
c_2=-\frac{65}{24}-\frac{21}{8}(c_1+h_1+k_1\ E_0)-k_2\ E_0 -h_2  \,  . 
\end{array}  \right. \monend  
These constraints  leave  as independent parameters :
$Pr$, $E_0$, $c_1$, $k_2$, $c_{\epsilon E}$, $c_{\varpi \rho}$, $c_{\varpi E}$ 
and the relaxation rates \quad $s_{XX}$, $s_{rx}$, $s_{\epsilon}$, $s_{\varpi}$ and $s_{{XXe}}$.
%
The free parameters are chosen to get a stable scheme by computing the roots of
the dispersion equation for several values of the wave vector ranging from $0$ to 
$ 2 \pi$ in magnitude and several directions with respect to the axis $x$ and $y$. 
Our approach is heuristic and nothing   {\it a priori} guaranties the $L^2$ stability. 

\bigskip \bigskip  \noindent {\bf \large 3) \quad Navier Stokes solver for  a compressible gas }  

\monitem 
This approach starts from the  
conservation equations of mass and momentum:
\moneq  \label{ns2d-eqs-conservation} 
 \left\{ \begin{array} {rl}  \displaystyle 
\frac{\partial}{\partial t}\rho + \frac{\partial}{\partial x} \rho v_x
+ \frac{\partial}{\partial y} \rho v_y &= 0  \\  \displaystyle 
 \frac{\partial}{\partial t}\rho v_x +  \frac{\partial}{\partial x} \rho v_x v_x
+ \frac{\partial}{\partial y}\rho  v_x v_y+\frac{\partial}{\partial x}P  
-\nu \, \triangle v_x -\zeta \, \frac{\partial}{\partial x} \big( {\rm div} \, v \big) &= 0 
 \\  ~   \vspace{-5 mm}   \\ \displaystyle
 \frac{\partial}{\partial t}\rho v_y + \frac{\partial}{\partial x} \rho v_x v_y
+ \frac{\partial}{\partial y}\rho v_y  v_y+\frac{\partial}{\partial y}P  
-\nu \, \triangle v_y -\zeta \, \frac{\partial}{\partial y} \big( {\rm div} \, v  \big) &= 0 \, . 
\end{array}  \right. \monend  
We add also the conservation of total energy. 
We assume the fluid is a perfect gas, then the  pressure $P$  is given 
according to 
\moneq  \label{ns2d-pression} 
P \,=\,  \rho \, R  \, T
\monend  
and the internal energy per unit mass $e$  is related to $T$ by 
\moneq  \label{ns2d-energie-interne} 
e  \,=\, {{R  T}\over{\gamma -1}} \, . 
\monend  
Then the evolution equation for the internal energy  takes the form 
\moneq  \label{ns2d-cons-energie} 
 \left\{ \begin{array} {l}  \displaystyle 
\frac{\partial}{\partial t}\rho e + \frac{\partial}{\partial x} \rho e v_x
+ \frac{\partial}{\partial y} \rho e v_y 
+P \Big[ \frac{\partial }{\partial x} u_x +\frac{\partial}{\partial y} u_y \Big] 
-\kappa \Big[\frac{\partial^2}{\partial x^2}T+\frac{\partial^2}{\partial y^2}T \Big] 
 \\  \displaystyle \qquad \quad  
-\nu \Big[ \Big( \frac{\partial}{\partial x}u_x-\frac{\partial}{\partial y}u_y \Big)^2
+ \Big( \frac{\partial}{\partial x}u_y+\frac{\partial}{\partial y}u_x \Big)^2 \Big]
-\zeta \, \Big(\frac{\partial}{\partial x}u_x+\frac{\partial}{\partial y}v_y \Big)^2 = 0  \,.
\end{array}  \right. \monend  

\monitem 
A linearized analysis gives the propagation and damping of the four hydrodynamic modes.
We deduce an algebraic expression for the     
sound velocity $c_s$, the relaxation $\nu_\tau$ of the transverse mode,  the relaxation $\nu_{\rm diff}$
 of the diffusive mode, and the damping  $\nu_{\rm acous}$ of the sound  modes: 
\moneq  \label{ns2d-relaxation} 
c_s \,= \, \sqrt{ \gamma R T} \,, \quad 
  \nu_\tau  \, = \,  \nu  \,, \quad 
\nu_{\rm diff} \, = \, \kappa \, \frac{\gamma -1}{R \gamma} \,,  \quad  
\nu_{\rm acous}   \, = \,   \frac{1}{2} \, (\nu + \zeta ) 
+ \frac{(\gamma -1)^2}{2 R \gamma} \,.     
\monend 
The non-linear terms allow to show that a uniform advection speed $\{V_x,V_y\}$ leads to 
phase shifts compatible with Galilean invariance. 

\monitem 
The model can be approximately simulated with simple 
finite difference expressions  for the space derivatives. 
We have developed a compressible Navier-Stokes solver for the numerical
resolution of the mathematical model (\ref{ns2d-eqs-conservation}) - 
(\ref{ns2d-cons-energie}). 
We use a cell vertex approach (with the nomenclature of Roache \cite{Ro72}). 
All the differential operators
are discretized with centered finite differences. The discrete evolution in times is obtained with an
elementary forward Euler first order explicit scheme. 
The  Dirichlet boundary conditions for velocity and temperature 
are implemented in a clear way by forcing the given value on the boundary node vertex.
For the adiabatic wall  where $ \, {{\partial T}\over{\partial n}} $ is null,  
a Neumann homogeneous boundary condition is enforced with mirror techniques decribed in the       
classical reference~\cite{Ro72}.

\bigskip \bigskip  \noindent {\bf \large 4) \quad  A simple test case }  

\monitem 
This test case has been studied in our contribution \cite{LD13}.     
The domain is a  $ \, N_x \times N_y \, $ rectangle  with periodic boundary conditions.
This test case is error-free as far as  boundary conditions are concerned. 
The initial  condition is a fluid at rest:  $ \, V_x = V_y = 0$.
The initial temperature $ \, T(x,y)= T_0 + \delta T_0 \, \cos { k \smb x} \, $
is associated with a wave number 
 $ \, k=2 \pi K / N_x$. 
Then density, pressure or energy are such that 
no acoustic wave is excited.
This is possible with the following conditions:     
\moneq  \label{ns2d-cns-initiale-densite} 
 \left\{ \begin{tabular} {rl}  
For D2Q9-D2Q5 : & \quad $  \rho=1  \,, $  \cr
For D2Q13 : &   \quad  $\rho = 1- 28 \ (T(x,y)-T_0)  \,,  $ \cr
For Navier Stokes :   & \quad $  P = R T_0 \,$, 
$ \, \rho E = {{P}\over{\gamma-1}} \,+\, {{\rho \, (V_x^2+V_y^2)}\over{2}} $.
\end{tabular}  \right. \monend  
One verifies that $T(x,y,t)$ relaxes exponentially in time.

\bigskip \bigskip  \newpage \noindent {\bf \large 5) \quad   Buoyancy }  

\smallskip 
   \centerline {\includegraphics[width=.61 \textwidth]   {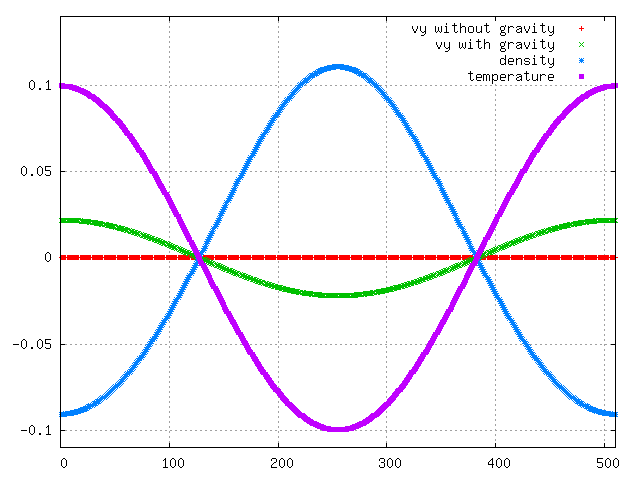} }  

\smallskip \noindent  {\bf Figure 2}. \quad 
Buoyancy  flow with the compressible Navier-Stokes solver.
The nonlinear exchanges between temperature and density are not affected by the gravity. 
\smallskip 

\monitem 
For the first scheme D2Q9-D2Q5 we simulate buoyancy
by adding  a vertical force ($V_y$) proportional to $T-T_0$. 
For the D2Q13 scheme and the direct approach of 
Navier-Stokes equations with finite differences, we
add vertical force in the $V_y$ momentum equation proportional to $\rho-\rho_0$.
Then the vertical speed increases approximately linearly with time 
and there is essentially no horizontal velocity.
With the  Navier-Stokes solver, 
we use a domain composed by  510 mesh points in width, and  periodic in height.
The temperature is periodic relative to the $x$ direction.     
With the D2Q13 lattice Boltzmann solver, we use the same domain as previously:
a domain of  510 meshes in width and periodic in height, with an (initial) 
temperature periodic in $x$ (see Fig.~3).

\centerline {\includegraphics[width=.60 \textwidth]    {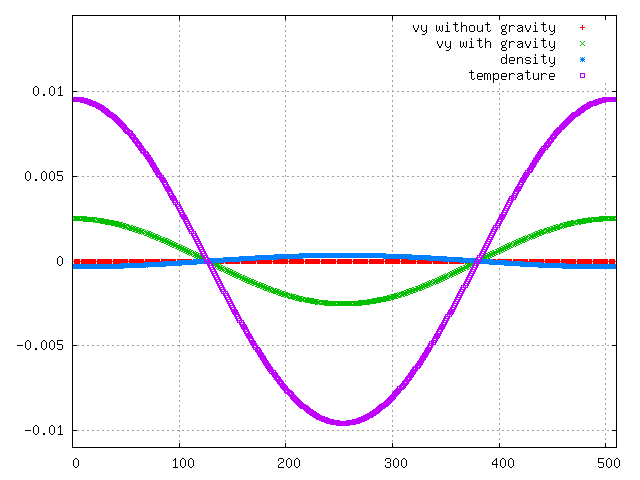} }  

\smallskip \noindent  {\bf Figure 3}. \quad 
Buoyancy  flow with the D2Q13 direct lattice Boltzmann solver. 
The nonlinear exchanges between temperature and density are not affected by the gravity. 
\smallskip 

\bigskip \bigskip  \noindent {\bf \large 6) \quad   Transverse waves }  

\monitem 
With the D2Q13 stencil, the non linear behaviour  for transverse waves    
is operating as follows. 
The initial conditions 
$ \, v_y(x,y,0) \equiv V_{y0} \, \cos{k  x} \, $
leads to density waves of wave vector $2 k$.
We modify the initial conditions by 
$ \, \rho(x,y,0) = \rho_0 + a\ v_y(x,y,0)^2 $. 
We measure the following global agregates relative to time: 
\moneq  \label{d2q13-mesure} 
\tilde{V}_y(t) \, = \, \sum_x v_y \, \cos ( k \ x)  \,, \quad 
\tilde{\rho}(t)  \, = \,  \sum_x \rho\,  \cos ( 2 \ k \ x )  \,, \quad 
\tilde{E}(t) \, = \,  \sum_x E \, \cos ( 2 \ k \ x)   \,.
\monend  
The typical result for D2Q13 is summarized in Fig.~4: 
the growth of $\tilde{E}$ in time  is proportional to $\nu \, k^2 \, V_y^2$.

\smallskip 
  \centerline {\includegraphics[width=.65 \textwidth]   {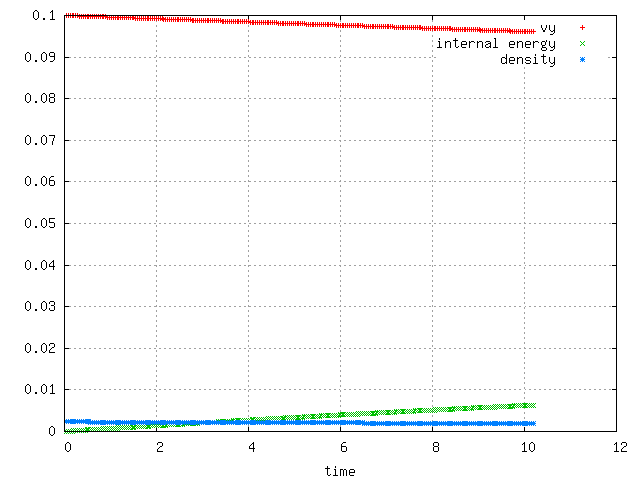} } 

\smallskip \noindent  {\bf Figure 4}. \quad 
Transverse waves with the D2Q13 direct lattice Boltzmann solver. 
The growth of $\tilde{E}$ in time  is proportional to $\nu \, k^2 \, V_y^2$.
\smallskip 

\monitem 
The interpretation of this evolution can be stated as follows.
The linearized equivalent equations at order 1 
with  space derivatives can be written in matrix form  
\moneq  \label{d2q13-equiv-edp-ordre-1} 
\begin{array}{|ccc|}
\partial_t &  \partial_r &  0 \cr
( {14\over13} + {1\over2}  V^2)\ \partial_r  & \partial_t & {1\over26}  \partial_r  \cr
0 &    {1\over2}  (39 + 13 c_1 + 2 E_0)\ \partial_r & \partial_t 
\end{array} \, =  \, 0 \, . \monend  
Without advective velocity, the diffusive mode is 
$ \, \,  ( 1, \, 0, \, -28 )^{\displaystyle \rm t} \, $
{\it id est}   $ \, E = -28 \, \rho $. 
With a transverse velocity $V$,  the diffusive mode is 
equal to  
$ \, \,  ( 1, \, 0, \, -28 - 13\, V^2  )^{\displaystyle \rm t} $,  {\it id est}\br 
 $  E = -28 \, \tilde{\rho} \, $ 
with a  density $\rho$    replaced by 
$ \, \tilde{\rho} \equiv  1 + {13\over28} V^2$. 

\bigskip \bigskip  \noindent {\bf \large 7) \quad   Couette flows  }  

\monitem 
Typical boundary conditions for Couette flows with the  D2Q9-D2Q5 scheme are stated as follows. 
For  $x=1$ and $x=N_x$ the velocity is known:    
$V_x=0$  and  $V_y$ is given. This type of boundary condition    
is classically achieved by a  ``bounce-back'' condition. 
Assuming a zero value for the mean temperature, for  $x=1$ the temperature is imposed:  $T=+\Delta T$  
and when  $x=N_x$ it has the opposite sign: $T=-\Delta T$ for $x=N_x$ .
This boundary condition is  achieved by an ``anti-bounce-back''
as proposed by  Ginzburg \cite {Gi05}.   
In consequence, the way we implement the boundary
conditions is not straightforward.
For the unit velocities with non-zero component parallel to the boundary,    
we use a bounce-back boundary condition. 
For the other velocities, an ``anti-bounce-back'' is implemented.
We consider for example: 
\moneq  \label{d2q13-cl-abb} 
f_1 + f_3  \, =  \, 2 \ ( p_\rho \ \rho + p_{XX} \ XX + p_E \ E)  
\monend  
with $E$ and $XX$ imposed and  $\rho$ estimated by extrapolation 
from values measured in the fluid.

 \centerline {\includegraphics[width=.55 \textwidth]   {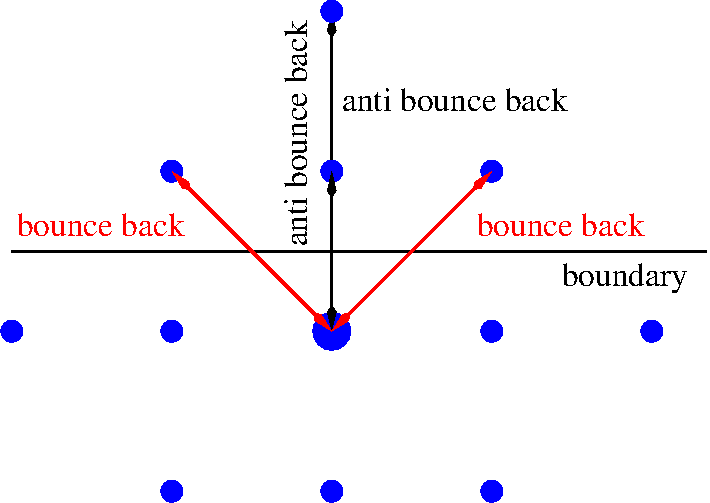} } 

\smallskip \noindent  {\bf Figure 5}. \quad 
Mixed   ``bounce-back''  and  ``anti-bounce-back'' boundary conditions for a flow    
simulated with the D2Q13 lattice Boltzmann scheme. 
In this case,  both velocity and temperature are imposed.
\smallskip 

\smallskip 
\centerline {\includegraphics[width=.65 \textwidth]   {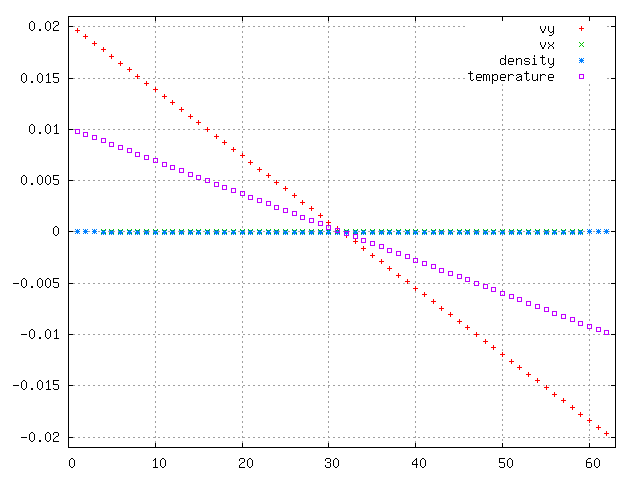}} 

\smallskip \noindent  {\bf Figure 6}. \quad 
Couette flow with the D2Q9-D2Q5 coupled scheme.
Result for the temperature field without motion of the lateral plates.  
\smallskip 

\monitem
With the coupled  D2Q9-D2Q5  scheme, the density is unchanged (see Fig.~6). For the
discretization of the compressible Navier Stokes equations, the pressure remains constant.  
Then, due to the equation of state (\ref{ns2d-pression}), 
the variation of density and temperature are coupled.
This effect is clearly visible in the Fig.~7 (direct Navier Stokes solver)      
With the D2Q13 lattice Boltzmann scheme, the internal energy $ \, e \, $ 
(proportional to the temperature) can be recovered 
thanks to the relation $ \, e \, = \,  E-\frac{13}{2}(V_x^2+V_y^2) \,$ 
as displayed in Fig.~8.                                                         

\smallskip 
  \centerline {\includegraphics[width=.65 \textwidth]   {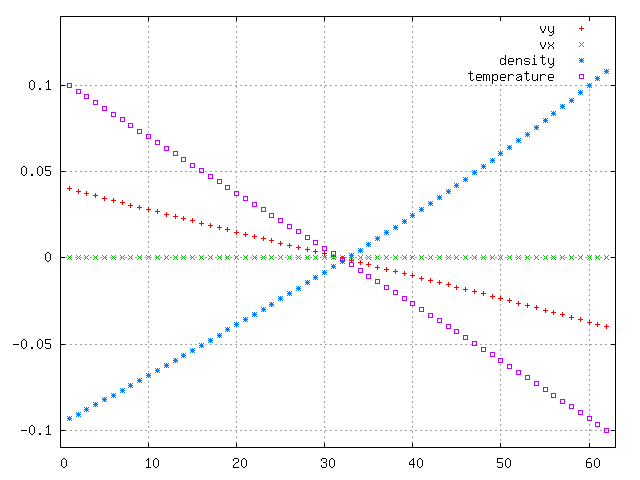}} 

\smallskip \noindent  {\bf Figure 7}. \quad 
Couette flow with a finite difference direct Navier-Stokes solver. 
Same result for the temperature field without motion of the plates.
\smallskip 

\smallskip 
\centerline {\includegraphics[width=.65 \textwidth]  {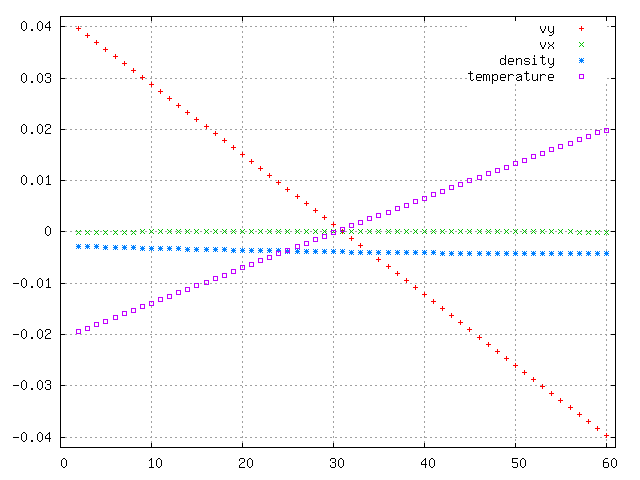} } 

\smallskip \noindent  {\bf Figure 8}. \quad 
Couette flow with  the D2Q13  lattice Boltzmann solver. 
The dashed curves show the ``temperature''  $\, T \equiv E-\frac{13}{2}(V_x^2+V_y^2) $.
There is no variation of the temperature 
when no gradient is imposed between the plates.
\smallskip 

\bigskip \bigskip  \newpage 
\noindent {\bf \large 8) \quad   Poiseuille flow    }  
 
\monitem
A Poiseuille flow is realized by adding an external term to take into account the gradient 
of pressure. 
Then a parabolic velocity profile is obtained as usual. 
Moreover, we add a Couette-type temperature profile between the lateral plates.  
In all our simulations, we do not observe any  
variation of the temperature  when no gradient is imposed between the plates. 
Moreover, when a  gradient of temperature is imposed, we observe a regular evolution  
of the temperature without destruction of the parabolic profile.

\smallskip 
  \centerline {\includegraphics[width=.65 \textwidth]  {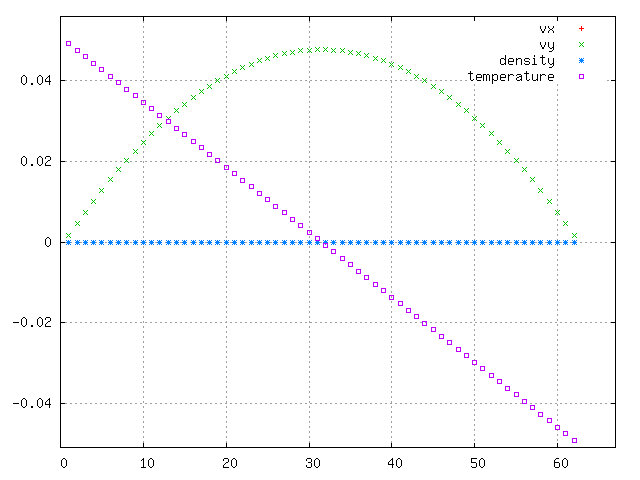} } 

\smallskip \noindent  {\bf Figure 9}. \quad 
Poiseuille  flow with the D2Q9-D2Q5 coupled scheme.
No variation of the temperature  when no gradient is imposed between the plates. 
\smallskip 

\smallskip 
  \centerline {\includegraphics[width=.65 \textwidth]  {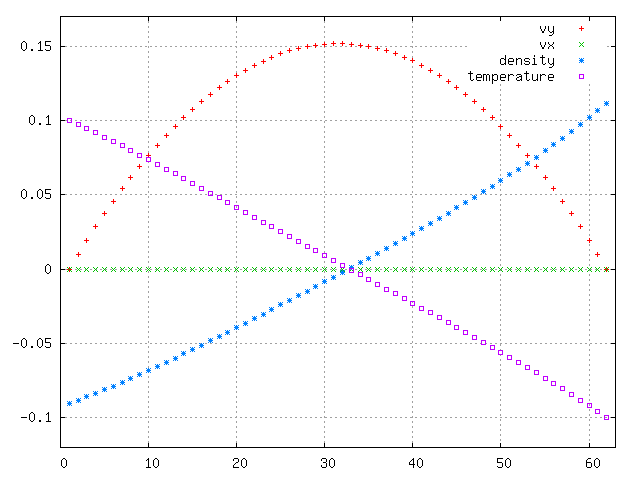} } 

\smallskip \noindent  {\bf Figure 10}. \quad 
Poiseuille  flow with a finite difference direct Navier-Stokes solver. 
No variation of the temperature   when no gradient is imposed between the plates.
\smallskip 

\smallskip 
\centerline {\includegraphics[width=.65 \textwidth]   {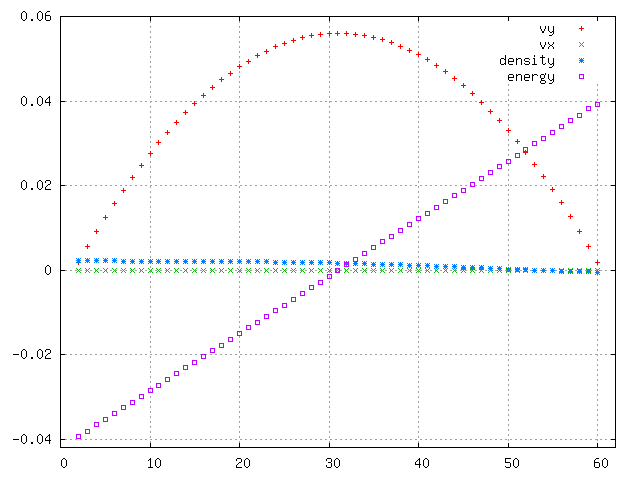} } 

\smallskip \noindent  {\bf Figure 11}. \quad 
Poiseuille  flow with  the D2Q13  lattice Boltzmann scheme. 
No variation of the temperature   when no gradient is imposed between the plates.
\smallskip 

\bigskip \bigskip  \noindent {\bf \large 9) \quad Test of an adiabatic boundary for the  D2Q13 scheme }  

\smallskip 
  \centerline {\includegraphics[width=.65 \textwidth]  {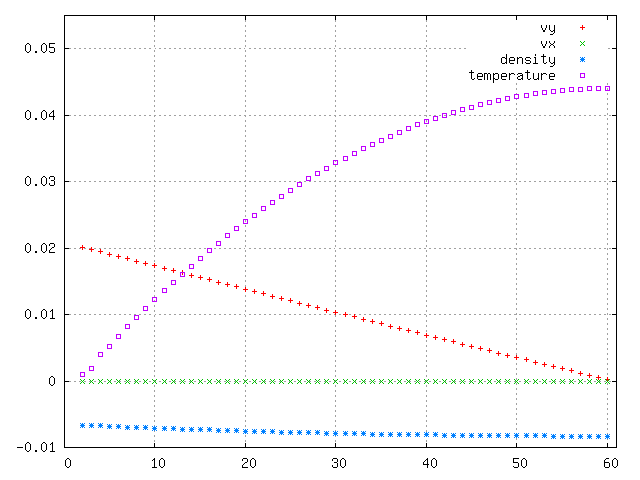} } 

\smallskip \noindent  {\bf Figure 12}. \quad 
Test of an adiabatic boundary with  the D2Q13  lattice Boltzmann solver. 
A uniform source of energy is applied.  
The $y$-velocity is null when  the left boundary is fixed.
The other fields are unchanged. 
\smallskip 

\monitem 
In order to implement correctly a null flux Neumann boundary condition
relative to the temperature,  we have tested our schemes for 
a uniform volumic source of energy. 
A homogeneous temperature given at the left boundary
and a homogeneous Neumann condition for temperature
 at the right  boundary. 
At the left  boundary the velocity is given as homogeneous or inhomogeneous:
$v_y = 0$ or  $v_y = 0.02$. 
We impose  a null velocity at the right  boundary. 
The solution is a ``semi-parabol'' and is correctly simulated as decribed in Fig.~12.   

\bigskip \bigskip  \noindent {\bf \large 10) \quad  Thermal test case of de Vahl Davis  }  

\monitem   
The de Vahl Davis \cite{dVD83} test has been described in the introduction.
We have used a $\,  187 \times  187 \, $ domain 
with a  Prandtl number equal to $0.71$ with the lattice Boltzmann simulations
and a  grid with $\,  256 \times  256 \, $ mesh points. 
For  the simple Navier-Stokes solver with finite differences, we have used 
$ \, 128 \times  128  , \, $
$ \, 196 \times  196   \, $ and 
$ \, 256 \times  256   \, $ mesh sizes. The results for the mean Nusselt number 
is respectively equal to 4.5099, 4.5154 and 4.5157 for a Rayleigh number equal to $10^5$. 
Our results are globally summarized in the following table. 
Simulations have been done on graphics card and implemented with Cuda. We do not have 
 precise comparisons of execution times between 
the two lattice Boltzmann models. We estimate the overhead to be roughly + 20 \% for  
D2Q13 compared to the coupled D2Q9-D2Q5.

\smallskip  \smallskip 
  \centerline {\includegraphics[width=1.3 \textwidth]  {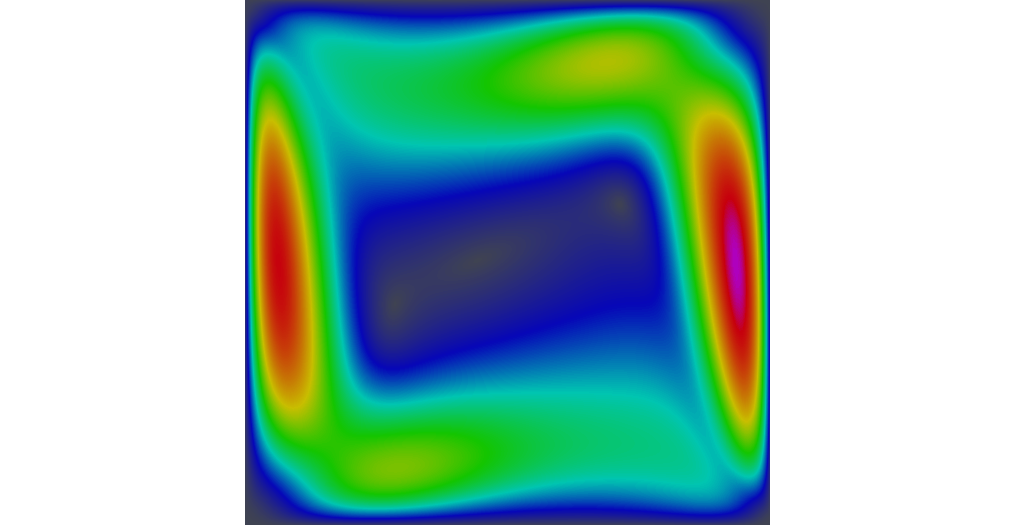} } 

\smallskip \noindent  {\bf Figure 13}. \quad 
De Vahl Davis thermal test case for natural convection
with  the D2Q13 direct lattice Boltzmann solver. 
Rayleigh number  = $10^5$.
Iso-velocity curves for the modulus  of the fluid speed.
The maximum velocity is 5.5\ $10^{-3}$.
\smallskip 

\newpage 
\smallskip \smallskip 

 \centerline { \begin{tabular}{|c|c|c|c|c|c|c|}    \hline 
 Rayleigh & de Vahl Davis &  Le Qu\'er\'e & Mezrhab &  D2Q9-D2Q5 & Navier Stokes & D2Q13   \\   \hline 
$10^5$          &  4.519   &      & 4.521 &  4.521  &  4.51  &   4.50   \\   \hline 
$10^6$          &  8.800 & 8.8252 & 8.824 &  8.828  &  8.88  &   8.73   \\   \hline  
\end{tabular} }  

\smallskip \noindent  {\bf Table 1}. \quad 
Comparison of Nusselt number integrated  in the whole cavity for two Rayleigh numbers. 

\smallskip 

\smallskip  \smallskip 
 \centerline {\includegraphics[width=1.3 \textwidth]    {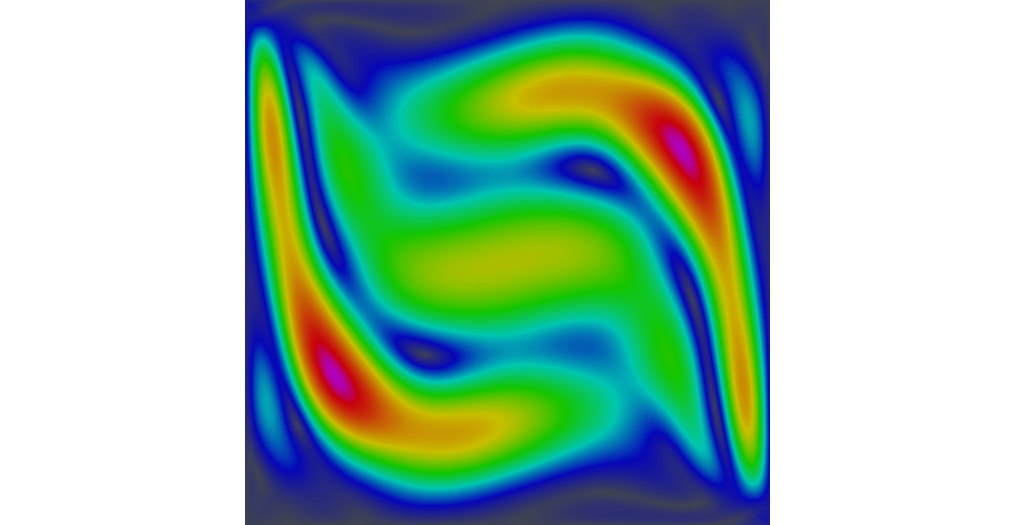}  } 

\smallskip \noindent  {\bf Figure 14}. \quad 
De Vahl Davis thermal test case for natural convection
with  the D2Q13 direct lattice Boltzmann solver. 
Rayleigh number = $10^5$.
Modulus of the asymmetry of the fluid speed.                                              
Modulus of the asymmetry of the Fluid Speed.
Curves for the departure from center symmetry. 
The maximum difference of velocity 
$|V(x,y)+V(x_0-x,y_0-y)|$ is 0.28\ $10^{-3}$ (5 \%).
\smallskip 

\monitem 
%
We have compared our results with those of de Vahl Davis \cite{dVD83},     
Le Qu\'er\'e \cite{LQ91}, Mezhrab  {\it et al.} \cite{MMJN10},  Wang {\it et al.} \cite{WWLL13}
with the coupled approach D2Q9-D2Q5 and our simple finite differences Navier-Stokes
solver. 
The results ``D2Q13'' obtained with a single particle distribution are correct but
must be considered as  preliminary compared to the other results.  
Inspection of the thermal and velocity fields obtained with the D2Q9-D2Q5 shows
that they are symmetric with respect to the center of the cavity. Similar inspection
for the fields obtained either with D2Q13 or the simple compressible Navier-Stokes
code used here show disymmetries that increase with the Rayleigh number. A detailed    
analysis of these asymmetries will be performed later and checked with data obtained   
with more sophisticated Navier-Stokes codes.

\bigskip \bigskip  \noindent {\bf \large Conclusion}  

\noindent

In this contribution, we have shown that coupled fluid and thermal flows that characterize  
natural convection can be simulated in two space dimensions with a single 
D2Q13 lattice Boltzmann scheme
by imposing   the conservation of mass, momentum and energy. 
We have tested our approach by a progressive  complexification of the test cases. 
Observe that strong compressible effects including the simulation of shock waves have not been    
considered in this contribution. 
The de Vahl Davis test case for natural convection 
gives encouraging results when we compare our result 
to previous ones obtained with a D2Q9-D2Q5 coupled approach or                            
with a direct simulation of the compressible Navier Stokes equations with finite
differences.  
Our results show that a lattice Boltzmann model  
with a single D2Q13 distribution that conserves mass, momentum and energy
gives results that compare better to direct Navier-Stokes simulations with finite 
differences than with simulations  obtained with the Boussinesq approximation.
Nevertheless, complementary studies are necessary to improve this method   
of simulation and confirm our results.

\bigskip \bigskip  \noindent {\bf \large Acknowledgments}   

 \noindent  
   The authors  thank the two anonymous referees for their very constructive remarks. 
Many thanks also to 
  the ``LaBS project'' (Lattice Boltzmann Solver,  www.labs-project.org), 
 funded by the French FUI8 research program, for supporting this contribution.


\medskip

\end{document}